\documentclass[12pt, a4paper]{article}

\usepackage{a4,amsmath,amssymb, amsthm, latexsym, color, graphicx,url}
\usepackage{tikz}
\usetikzlibrary{arrows.meta}
\usetikzlibrary{decorations.markings}
\usetikzlibrary{calc}

\usepackage{fullpage}
\usepackage{enumitem}
\usepackage{xcolor}
\newtheorem{theorem}{Theorem}[section]
\newtheorem{proposition}[theorem]{Proposition}
\newtheorem{lemma}[theorem]{Lemma}
\newtheorem{corollary}[theorem]{Corollary}

\theoremstyle{definition}
\newtheorem{example}[theorem]{Example}
\newtheorem{definition}[theorem]{Definition}
\newtheorem{remark}[theorem]{Remark}

\newtheorem{openproblem}[theorem]{Open Problem}

\title{New constructions for disjoint partial difference families and external partial difference families}
\author{Sophie Huczynska and Laura Johnson \thanks{email: sh70@st-andrews.ac.uk,  lj68@st-andrews.ac.uk}}
\date{School of Mathematics \& Statistics, University of St Andrews, St Andrews, KY16 9SS, Scotland, UK }
\begin{document}
\maketitle

\begin{abstract}
Recently, new combinatorial structures called disjoint partial difference families (DPDFs) and external partial difference families (EPDFs) were introduced, which simultaneously generalize partial difference sets, disjoint difference families and external difference families, and have applications in information security.  So far, all known construction methods have used cyclotomy in finite fields. We present the first non-cyclotomic infinite families of DPDFs which are also EPDFs, in structures other than finite fields (in particular cyclic groups and non-abelian groups). As well as direct constructions, we present an approach to constructing DPDFs/EPDFs using relative difference sets (RDSs); as part of this, we demonstrate how the well-known RDS result of Bose extends to a very natural construction for DPDFs and EPDFs.
\end{abstract}

\section{Introduction}
Difference sets and difference families are well-studied combinatorial objects dating back to the 1930s; difference families are useful for constructing balanced incomplete block designs (BIBDs) (\cite{ColDin}, \cite{Wil}).  Disjoint difference families (DDFs) have recently received attention \cite{Bur}, with applications to design theory and information security.  In the early 2000s, motivated by applications in cryptography, external difference families (EDFs) were introduced (\cite{OgaKurStiSai}, \cite{PatSti}).  In a recent paper \cite{HucJoh}, partial analogues of DDFs and EDFs were introduced; these are called Disjoint Partial Difference Families (DPDFs) and External Partial Difference Families (EPDFs). A $(v,s,k,\lambda,\mu)$-DPDF (respectively, EPDF) is a set $S$ of $s$ disjoint $k$-subsets of an order-$v$ group, such that the multiset union of internal (respectively, external) differences of the sets in $S$ comprises $\lambda$ copies of each non-identity element in $S$, and $\mu$ copies of each non-identity element not in $S$.  These also generalize the concept of a partial difference set (PDS) (see \cite{Ma84}, \cite{Ma}) and have applications in information security. In \cite{HucJoh}, construction methods were given for DPDFs and EPDFs in $GF(q)$ (where $q$ is a prime power) using cyclotomic techniques.    Cyclotomy has long been used to produce traditional difference families, beginning with the work of \cite{Wil}.  This paper takes the first step in going beyond cyclotomy to present a range of other construction methods in structures other than finite fields.

It is of particular interest to construct families of sets which are simultaneously DPDFs and EPDFs.  It is shown in \cite{HucJoh} that such families must partition a partial difference set, which is regular if it is proper.  As well as the natural theoretical appeal of such examples, DPDFs which partition a regular PDS correspond to a two-class association scheme which means they can be used to obtain partially balanced incomplete block designs (PBIBDs) (\cite{Joh}).

In \cite{HucJoh}, DPDFs/EPDFs are obtained by partitioning cyclotomic PDSs in finite fields; for fields of prime order, these partition the quadratic residues (or non-residues) modulo $p$.  In this paper, we address the following questions.
\begin{itemize}
    \item Can DPDF/EPDF constructions be obtained in abelian groups other than $(GF(q),+)$; particularly for cyclic groups $\mathbb{Z}_v$ where $v$ is not prime?
    \item Can DPDF/EPDF constructions be obtained in non-abelian groups?   
\end{itemize}
We provide a range of constructions answering both of these questions in the affirmative.

It is known (\cite{Ma84}, \cite{Ma}) that if $D$ is a regular PDS in $\mathbb{Z}_v$, then there are just two possibilities: $D$ or its complement is the set of non-zero squares (equivalently, non-squares) in $GF(v)$ with $v \equiv 1 \mod 4$ prime, or else $D \cup \{0\}$ or $G \setminus D$ is a subgroup of $G$.  Constructions of DPDFs and EPDFs partitioning the former type of PDS was addressed in \cite{HucJoh}; in this paper we address the latter situation (although not limited to the cyclic group setting).  In non-abelian groups, while the definitions of these difference family-type structures remain valid, very little is known.  There are just a few non-abelian EDFs in the literature (see for example \cite{HucJefNep}), and prior to this paper there were no known constructions for non-abelian DPDFs or EPDFs, so the non-abelian constructions presented here are significant.

We first present explicit constructions in cyclic groups. We next develop constructions in general finite groups based on relative difference sets (RDSs), in which the subgroup not present in the union of the sets of the DPDF/EPDF is precisely the forbidden subgroup for the RDS.  In particular, we show how the classic result of Bose \cite{Bos} which originally constructed relative difference sets using finite geometry, very naturally extends to a DPDF/EPDF construction in cyclic groups.  We obtain a framework for using RDSs for DPDF/EPDF constructions which encompasses this example and generates many others.  Finally, we briefly present DPDFs and EPDFs in cyclic groups which demonstrate that not all DPDFs must be EPDFs, and vice versa.

Any $(v,s,k,\lambda,0)$-DPDF which partitions $G\setminus H$ (for some subgroup $H$ of $G$) gives an instance of a \emph{relative difference family}.  Relative difference families were introduced in \cite{Bur98} and subsequently explored in various further papers (eg \cite{Bur99}, \cite{Mom}).  These have mostly been studied in abelian groups, and are closely related to the concept of \emph{group divisible designs}. EPDFs also give examples of \emph{bounded external difference families} (see \cite{PatSti}).

\section{Background}
Throughout what follows, we let $G$ be a group, written additively unless otherwise stated, and let $G^*$ denote $G \setminus \{0\}$.  

For a subset $D$ of $G$, we define the multiset 
$$\Delta(D)=\{ x-y: x \neq y \in D \}$$
and for sets $D_1, D_2 \subseteq G$, we define the multiset 
$$\Delta(D_1, D_2)=\{ x-y: x \in D_1, y \in D_2 \}.$$  
(In multiplicative notation these are $\Delta(D)=\{ xy^{-1}: x \neq y \in D \}$ and $\Delta(D_1, D_2)=\{ xy^{-1}: x \in D_1, y \in D_2 \}$ respectively.)\\
For a family $A=\{A_1, \ldots, A_s\}$ of disjoint subsets of $G$, we define 
$$ {\rm Int}(A)=\bigcup_{i=1}^s \Delta(A_i)$$
and
$$ {\rm Ext}(A)=\bigcup_{1 \leq i \neq j \leq s} \Delta(A_i,A_j).$$
For $g \in G$ and $S \subseteq G$, we denote the translate $g+S=\{g+s: s \in S\}$ (multiplicatively, $gS=\{gs: s \in S\}$).

We begin with a summary of relevant definitions (see \cite{ColDin} and \cite{HucJoh}).
\begin{definition}
Let $G$ be a group of order $v$.
\begin{itemize}
\item[(i)] A $(v,k,\lambda, \mu)$-partial difference set (PDS) is a $k$-subset $D$ of $G$ with the property that the multiset of differences $\Delta(D)$ comprises each non-identity element of $D$ precisely $\lambda$ times, and each non-identity element of $G \setminus D$ $\mu$ times. If $\lambda=\mu$ then $D$ is simply called a $(v,k,\lambda)$-difference set (DS); otherwise the PDS is said to be proper.
\item[(ii)] A $(v,s,k,\lambda)$-disjoint difference family (DDF) is a collection of disjoint $k$-subsets $S^{\prime}=\{A_1,\ldots,A_s\}$ of $G$ with the property that ${\rm Int}(S^{\prime})$ comprises each non-identity element of $G$ precisely $\lambda$ times.   If the disjointness condition is relaxed we obtain a difference family (DF).
\item[(iii)] A $(v,s,k,\lambda,\mu)$-disjoint partial difference family (DPDF) is a collection of disjoint $k$-subsets $S^{\prime}=\{A_1,\ldots,A_s\}$ of $G^*$ with the property that ${\rm Int}(S^{\prime})$ comprises each non-identity element of $S=\cup_{i=1}^s A_i$ precisely $\lambda$ times, and each non-identity element of $G \setminus S$ $\mu$ times.  If $\lambda=\mu$ then $S^{\prime}$ is a $(v,s,k,\lambda)$-DDF.
\item[(iv)] A $(v,s,k,\lambda)$-external difference family (EDF) is a collection of disjoint $k$-subsets $S^{\prime}=\{A_1,\ldots,A_s\}$ of $G^*$ with the property that ${\rm Ext}(S^{\prime})$ comprises each non-identity element of $G$ precisely $\lambda$ times.  An EDF which partitions $G^*$ is called \emph{near-complete}.  
\item[(iii)] A $(v,s,k,\lambda,\mu)$-external partial difference family (EPDF) is a collection of disjoint $k$-subsets $S^{\prime}=\{A_1,\ldots,A_s\}$ of $G^*$ with the property that ${\rm Ext}(S^{\prime})$ comprises each non-identity element of $S=\cup_{i=1}^s A_i$ precisely $\lambda$ times, and each non-identity element of $G \setminus S$ $\mu$ times.  If $\lambda=\mu$ then $S^{\prime}$ is a $(v,s,k,\lambda)$-EDF.
\end{itemize}
\end{definition}

\begin{lemma}
\begin{itemize}
    \item[(i)] If $S^{\prime}$ is a $(v,s,k,\lambda_1,\mu_1)$-DPDF then
    \begin{equation}\label{DPDFeqn}
        sk(k-1)=\lambda_1 sk + \mu_1 (v-1-sk)
    \end{equation}
    \item[(ii)] If $S^{\prime}$ is a $(v,s,k,\lambda_2,\mu_2)$-EPDF then  
    \begin{equation}\label{EPDFeqn}
    s(s-1)k^2= \lambda_2 sk + \mu_2 (v-1-sk)
    \end{equation}
\end{itemize}
\end{lemma}
\begin{proof}
This is immediate upon double-counting the elements of ${\rm Int}(S^{\prime})$ and ${\rm Ext}(S^{\prime})$.
\end{proof}

We will also need the definition of a relative difference set, and its generalization, the divisible difference set.  For a comprehensive survey article on these structures, see \cite{Pot}.
\begin{definition}\label{RDSdefinition}
Let $G$ be a group of order $mn$ and let $H$ be a normal subgroup of $G$ of order $n$. A $k$-subset $R$ of $G$ is an $(m,n,k,\lambda)$-relative difference set (RDS) in $G$ relative to $H$ if the multiset $\Delta(R)$ comprises each element in $G \setminus H$ exactly $\lambda$ times, and each non-identity element in $H$ exactly $0$ times.  If $n=1$ then $R$ is a difference set.
\end{definition}

A counting argument shows that, for an $(m,n,k,\lambda)$-RDS, we have the relation $k(k-1)=(mn-n)\lambda$.

\begin{definition}
Let $G$ be a group of order $mn$ and let $H$ be a normal subgroup of $G$ of order $n$.  A $k$-subset $D$ of $G$ is an $(m,n,k, \lambda,\mu)$-divisible difference set (DDS) relative to $H$ if the multiset $\Delta(D)$ comprises each non-identity element of $H$ exactly $\lambda$ times, and each element of $G \setminus H$ exactly $\mu$ times.  If $\mu=0$ then $D$ is a relative difference set.
\end{definition}
In general, more is known about RDSs than about DDSs.

\begin{remark}
In Definition \ref{RDSdefinition}, it is possible to relax the condition that $H$ is a normal subgroup.  An example of an RDS in $A_5$ relative to a subgroup $H$ of order $2$ is presented in \cite{CheLi}, which satisfies all conditions of Definition \ref{RDSdefinition}, except for the requirement that $H$ is normal (this would be impossible since $A_5$ is a simple group).
\end{remark}

\subsection{DPDFs and EPDFs partitioning PDSs}

In this section, we will explore the special properties of DPDFs and EPDFs which partition PDSs.

The following key result was proved in \cite{HucJoh}:
\begin{theorem}\label{PDSpartition}
Let $G$ be a group of order $v$.  Let $A=\{A_1, \ldots, A_s\}$ be a family of disjoint subsets of $G^*$, each of size $k$.  Then any two of the following conditions implies the third:
\begin{itemize}
    \item[(i)] $A$ partitions a $(v,k,\lambda,\mu)$-PDS in $G$;
    \item[(ii)] $A$ is a $(v,s,k,\lambda_1, \mu_1)$-DPDF in $G$;
    \item[(iii)] $A$ is a $(v,s,k, \lambda-\lambda_1, \mu-\mu_1)$-EPDF in $G$.
\end{itemize}
\end{theorem}
Moreover, if the PDS in Theorem \ref{PDSpartition} is proper, then it is regular.

As mentioned in the introduction, results have been obtained (\cite{Ma84}, \cite{Ma}) which significantly restrict the possibilities for regular PDSs. For cyclic groups, the following holds:

\begin{theorem}
Let $\mathbb{Z}_v$ be the cyclic group of order $v$. Let $S^{\prime}$ be a $(v,s,k,\lambda_1,\mu_1)$-DPDF and a $(v,s,k,\lambda_2,\mu_2)$-EPDF in $\mathbb{Z}_v$ which partitions a proper PDS.  Then 
\begin{itemize}
    \item[(i)] if $v$ is a prime and $v \equiv 3 \mod 4$ then no such $S^{\prime}$ exists;
    \item[(ii)] if $v$ is a prime and $v \equiv 1 \mod 4$ then $S^{\prime}$ partitions the set of non-zero quadratic residues or the non-residues modulo $v$;
    \item[(iii)] if $v$ is a composite number then $S^{\prime}$ partitions a proper non-trivial subgroup $H$ of $\mathbb{Z}_v$ or its complement $\mathbb{Z}_v \setminus H$.
\end{itemize}
\end{theorem}
\begin{proof}
It is known (\cite{Ma84}, \cite{Ma}) that if $G$ is a cyclic group of order $v$ and $D$ is a regular PDS in $G$ then either $v$ is an odd prime such that $v \equiv 1 \mod 4$ and $D$ is the set of quadratic residues (or non-residues) modulo $v$; or $D \cup \{0\}$ or $G \setminus D$ is a subgroup of $G$.  
\end{proof}

Examples of DPDFs and EPDFs partitioning a PDS of each type are given below:
\begin{example}\label{DPDFexample}
\begin{itemize}
\item [(i)]  Let $G=\mathbb{Z}_{13}$; the sets
$$ \{1,3,9\}, \{4,10,12\}$$
form a (13, 2, 3, 0, 2)-DPDF and a (13, 2, 3, 2, 1)-EPDF which partition the quadratic residues mod $13$ (see \cite{HucJoh}).
\item[(ii)] Let $G=\mathbb{Z}_{16}$ and $H=\{0,4,8,12\} \leq G$; the sets
$$\{1,9\},\{5,13\},\{2,14\},\{6,10\},\{3,15\},\{7,11\} $$
form a $(16,6,2,0,4)$-DPDF and $(16,6,2,14,8)$-EPDF which partition $G \setminus H$ (see Theorem \ref{cosetpartition}).
\end{itemize}
\end{example}

The following basic PDS result is useful (see \cite{Ma}):
\begin{lemma}\label{PDSlemma}
Let $G$ be a group of order $mn$ with identity $e$ and subgroup $H$ of order $n$.
\begin{itemize}
    \item[(i)] The sets $H$, $H \setminus \{e\}$, $G \setminus H$ and $(G \setminus H) \cup \{e\}$ are partial difference sets, with $H \setminus \{e\}$ and $G \setminus H$ being regular. 
    \item[(ii)] $H \setminus \{e\}$ is an $(mn, n-1, n-2, 0)$-PDS. 
    \item[(iii)] $G \setminus H$ is an $(mn,mn-n, mn-2n,mn-n)$-PDS.
 \end{itemize}
 \end{lemma}
    
In the case when a regular PDS $D$ is a subgroup with the identity removed, then we can characterize any DPDF or EPDF which partitions $D$, as follows.

\begin{theorem}
Let $G$ be a group of order $mn$ and $H$ be a subgroup of $G$ of order $n$.  Let $D=H \setminus \{0\}$.
\begin{itemize}
    \item[(i)] If $S^{\prime}$ is an $(mn,s,k,\lambda,\mu)$-DPDF (respectively, EPDF) partitioning $D$, then $\mu=0$ and $S^{\prime}$ is a near-complete $(n,s,k,\lambda)$-DDF (respectively, EDF) in the group $H$. 
    \item[(ii)] Each near-complete $(n,s,k,\lambda)$-DDF (respectively, EDF) in the group $H$ corresponds to an $(mn,s,k,\lambda,0)$-DPDF (respectively, EPDF) in $G$ partitioning $D$.
\end{itemize}
\end{theorem}
\begin{proof}
For (i), we establish the DDF case; the EDF case then follows by Theorem \ref{PDSpartition}.  Let $S^{\prime}$ be an $(mn,s,k,\lambda,\mu)$-DPDF partitioning $D$.  By definition, ${\rm Int}(S^{\prime})$ must comprise every element of $D$ (i.e. every non-identity element of $H$) $\lambda$ times and every non-identity element of $G \setminus D$ (i.e. $G \setminus H$) $\mu$ times.  Since $S=D \subseteq H$, all elements of ${\rm Int}(S^{\prime})$ lie in $H^*$, and so $\mu=0$.  Thus  ${\rm Int}(S^{\prime})$ comprises $\lambda$ copies of the non-identity elements of $H$, and $S^{\prime}$ partitions $H^*$, so $S^{\prime}$ is a near-complete DDF in $H$.  Correspondingly $S^{\prime}$ is also a near-complete EDF in $H$.\\
Part (ii) is clear, using the natural embedding of $H$ into $G$.
\end{proof}

\begin{example}
It can be verified that $\{1,4\}, \{2,3\}$ form a $(5,2,2,1)$-DDF and $(5,2,2,2)$-EDF in $\mathbb{Z}_5$. \\ 
The group $\mathbb{Z}_{10}$ contains the subgroup $H=\{0,2,4,6,8\} \cong \mathbb{Z}_5$, via embedding $f: \mathbb{Z}_5 \rightarrow \mathbb{Z}_{10}$, $x \mapsto 2x$.
Then $\{2,8\}, \{4,6\}$ is a $(10,2,2,1,0)$-DPDF and $(10,2,2,2,0)$-EPDF in $\mathbb{Z}_{10}$.
\end{example}

Since EDFs have been well-studied elsewhere (see \cite{PatSti} and references therein), we therefore focus on the situation when the PDS $D$ is the complement of a subgroup of $G$.

The next result guarantees that there exists a  DPDF/EPDF of this type, in any group $G$ containing a normal subgroup $H$.

\begin{theorem}\label{DPDFcosets}
Let $G$ be a group of order $mn$ and $H$ a normal subgroup of $G$ of order $n$.  Then the set of cosets of $H$ in $G$, excepting $H$ itself, forms an $(mn,m-1,n,0,mn-n)$-DPDF and an $(mn,m-1,n,mn-2n,0)$-EPDF.
\end{theorem}
\begin{proof}
By Lemma \ref{PDSlemma}, $G \setminus H$ is an $(mn,mn-n,mn-2n,mn-n)$-PDS.  The cosets of $H$ in $G$, other than $H$ itself, partition $G \setminus H$ and for each, its internal difference multiset comprises $n$ copies of $H^*$ and $0$ copies of $G \setminus H$.  So these cosets from an $(mn,m-1,n,0,mn-n)$-DPDF and consequently an $(mn,m-1,n,mn-2n,0)$-EPDF by Theorem \ref{PDSpartition}.
\end{proof}

We end this section with a result about the possible parameters for DPDFs/EPDFs which partition the complement of a subgroup.  We first need a technical lemma which in fact applies more widely to any DPDF or EPDF; note however it will never apply to any of the DPDFs or EPDFs from \cite{HucJoh} which partition a cyclotomic class.

\begin{lemma}\label{lemmacoprimeness}
If $\gcd(sk,v-1)=1$ then
\begin{itemize}
    \item[(i)] For a $(v,s,k,\lambda_1,\mu_1)$-DPDF, either $\mu_1=0$ or $\mu_1=sk$.
    \item[(ii)] For a $(v,s,k,\lambda_2,\mu_2)$-EPDF, either $\mu_2=0$ or $\mu_2=sk$.
\end{itemize}
\end{lemma}
\begin{proof}
For (i), we use Equation (\ref{DPDFeqn}); rearranging we see that
$$ sk(k-1)=sk(\lambda_1-\mu_1)+ \mu_1(v-1).$$
Hence $sk|\mu_1 (v-1)$; since $sk$ and $v-1$ are coprime, $sk|\mu_1$, but note $\mu_1 \leq sk$ since there are $sk$ elements in the sets of $S^{\prime}$.  Hence $\mu_1=0$ or $\mu_1=sk$.  Part (ii) follows from a similar rearrangement of Equation (\ref{EPDFeqn}). 
\end{proof}

\begin{theorem}\label{GsetminusH}
Let $G$ be a group of order $v=mn$.
Suppose $S^{\prime}$ is a $(v,s,k,\lambda_1,\mu_1)$-DPDF and a $(v,s,k,\lambda_2,\mu_2)$-EPDF that partitions $G \setminus H$ where $H \leq G$ has order $n$.  Then 
\begin{itemize}
    \item[(i)] $n \mid \mu_1$ and $n \mid \mu_2$.
    \item[(ii)] If $\gcd(mn-n,mn-1)=1$ then $S^{\prime}$ is one of the following:
    \begin{itemize}
    \item[(a)] an $(mn,s,k,k-1,0)$-DPDF and an $(mn,s,k,mn-2n-k+1,mn-n)$-EPDF;
    \item[(b)] an $(mn,s,k,k-n,mn-n)$-DPDF and an $(mn,s,k,mn-n-k,0)$-EPDF.
    \end{itemize}
\item[(iii)] If $n=2$ then $S^{\prime}$ is one of the following:
\begin{itemize}
\item[(a)] an $(2m,s,k,k-1,0)$-DPDF and an $(2m,s,k,2m-3-k,2m-2)$-EPDF;
\item[(b)] an $(2m,s,k,k-2,2m-2)$-DPDF and an $(2m,s,k,2m-2-k,0)$-EPDF.
\end{itemize}
\end{itemize}
\end{theorem}
\begin{proof}
For (i), using the fact that $sk=mn-n$ and $v-1=mn-1$ in Equation (\ref{DPDFeqn}), we have that
$$n(m-1)(k-1)=\lambda_1 n(m-1) + \mu_1 (n-1).$$
So $n|\mu_1(n-1)$, but since $\gcd(n-1,n)=1$, we must have $n|\mu_1$.  We may apply a similar argument using Equation (\ref{EPDFeqn}) to see $n|\mu_2$. Part (ii) is an application of Lemma \ref{lemmacoprimeness}.  For part (iii), $\gcd(sk,v-1)=\gcd(mn-n,mn-1)$ and any common divisor of $mn-n$ and $mn-1$ divides $mn-1-(mn-n)=n-1$, i.e. these quantities are coprime when $n=2$.
\end{proof}

Constructions producing DPDFs/EPDFs of type(a) in Theorem \ref{GsetminusH}(ii)/(iii) are presented in Section 4 using RDSs, while Section 3 includes constructions giving DPDFs/EPDFs of type (b) in cyclic groups.\\

\section{Cyclic DPDFs/EPDFs}

We present constructions for infinite families of DPDFs/EPDFs such that $G$ is a cyclic group, $H$ is a subgroup of $G$ and the DPDF/EPDF partitions $G \setminus H$. 

In this section, we will use the group ring notation whereby $\lambda S$ indicates the multiset comprising $\lambda$ copies of a set $S$.  (In general we shall avoid this notation in the rest of the paper, to avoid confusion with multiplicative translates of a set.)

We first introduce a family of subsets $S_i$ of $\mathbb{Z}_{2m}$ ($m>3$ odd) which have the useful property that $\Delta(S_i)$ and $\Delta(S_i,S_j)$ consist entirely of unions of $S_k$'s and copies of $H \setminus \{0\}$.  

\begin{proposition}\label{Siprop}
Let $m>3$ be an odd integer and let $G=\mathbb{Z}_{2m}$.  Let $H=\{0,m\}$ be the order-$2$ subgroup of $G$.\\
For $1 \leq i \leq 2m-1$, define
$$ S_i=\{i,m-i,m+i,2m-i\} \subseteq \mathbb{Z}_{2m}.$$
Then
\begin{itemize}
    \item[(i)] For $1 \leq i \leq 2m-1$,
    $$S_i=S_{m-i}=S_{m+i}=S_{2m-i}.$$\\
    In particular $S_1, S_2,\ldots, S_{\frac{m-1}{2}}$ comprise all the distinct $S_j$ in $G$ ($1 \leq j \leq 2m-1$).
    \item[(ii)] Each $|S_i|=4$ and $\{S_1,\ldots, S_{\frac{m-1}{2}} \}$ partition $G \setminus H$.
    \item[(iii)] $$\Delta(S_i)=4\{m\} \cup 2 S_{2i}.$$
    \item[(iv)] $$\Delta(S_i,S_j)= 2S_{i-j} \cup 2 S_{i+j}.$$
\end{itemize}
\end{proposition}
\begin{proof}
Part (i) is immediate from the definition of $S_i$.  For part (ii), the fact that $m$ is odd guarantees that all $4$ elements are distinct.  It is clear that as $i$ runs through $1,\ldots, \frac{m-1}{2}$, the sets $S_i$ account for all non-identity elements of $G$ other than $m$.  \\
For part (iii), write $S_i=A_i \cup B_i$ where $1 \leq i \leq \frac{m-1}{2}$, $A_i=\{i,m+i\}$ and $B_i=\{m-i,2m-i\}$. Then $\Delta(S_i)=\Delta(A_i) \cup \Delta(A_i,B_i) \cup \Delta(B_i,A_i) \cup \Delta(B_i)$ where $\Delta(A_i)=\Delta(B_i)=\{m,m\}$, $\Delta(A_i,B_i)=\{2i,2i,m+2i,m+2i\}=2A_{2i}$ and $\Delta(B_i,A_i)=2B_{2i}$.\\
For part (iv), $\Delta(S_i,S_j)=\Delta(A_i,A_j) \cup \Delta(A_i,B_j) \cup \Delta(B_i,A_j) \cup \Delta(B_i,B_j)$ where
$\Delta(A_i,A_j)=2A_{i-j}$, $\Delta(A_i,B_j)=2A_{i+j}$, $\Delta(B_i,A_j)=2B_{i+j}$ and $\Delta(B_i,B_j)=2B_{i-j}$.  Hence $\Delta(S_i,S_j)=2S_{i-j} \cup 2S_{i+j}$.
\end{proof}

\begin{theorem}\label{directcyclic}
Let $m>3$ be an odd number.  Let $G=\mathbb{Z}_{2m}$ and $H=\{0,m\} \leq G$.\\ 
The family of sets $S^{\prime}=\{S_i: 1 \leq i \leq \frac{m-1}{2} \}$ in $G$ given by
$$ S_i=\{i,m-i,m+i,2m-i \}$$
is a $(2m, \frac{m-1}{2}, 4,2,2m-2)$-DPDF and a $(2m, \frac{m-1}{2}, 4, 2m-6, 0)$-EPDF which partitions $G \setminus H$.
\end{theorem}
\begin{proof}
By Proposition \ref{Siprop}, the family $S^{\prime}$ partitions $S=G \setminus H$.
We consider ${\rm Int}(S^{\prime})= \cup_{i=1}^{\frac{m-1}{2}}\Delta(S_i)$.  Again by Proposition \ref{Siprop}, the multiset $\Delta(S_i)$ comprises 4 copies of $\{m\}$ and 2 copies of $S_{2i}$.  So $${\rm Int}(S^{\prime})= 4\left(\frac{m-1}{2}\right) \{m\} \cup 2(S_2 \cup S_4 \cup \cdots S_{m-3} \cup S_{m-1}).$$
If $m \equiv 1 \mod 4$, 
$$ S_2 \cup S_4 \cup \cdots \cup S_{m-1}=S_2 \cup S_4 \cup \cdots S_{\frac{m-1}{2}} \cup S_{m-\frac{m-3}{2}} \cup \cdots \cup S_{m-1}$$
while if $m \equiv 3 \mod 4$, 
$$ S_2 \cup S_4 \cup \cdots \cup S_{m-1}=S_2 \cup S_4 \cup \cdots S_{\frac{m-3}{2}} \cup S_{m-\frac{m-1}{2}} \cup \cdots \cup S_{m-1}.$$
In either case this union equals $S_1 \cup S_2 \cup S_3 \cup \cdots S_{\frac{m-1}{2}}=S$. Hence ${\rm Int}(S^{\prime})$ comprises $2m-2$ copies of $\{m\}$ and $2$ copies of $G \setminus H$, so is a $(2m,\frac{m-1}{2}, 4, 2, 2m-2)$-DPDF.  Since $G \setminus H$ is a $(2m,2m-2,2m-4,2m-2)$-PDS, $S^{\prime}$ is a $(2m,\frac{m-1}{2}, 4, 2m-6, 0)$-EPDF.
\end{proof}

For $m=3$, Theorem \ref{directcyclic} can still be used but, rather than constructing a family of sets, it yields just one set $\{1,2,4,5\}$ which is an RDS.  

\begin{example}\label{directcyclicex}
\begin{itemize}
    \item[(i)] Applying Theorem \ref{directcyclic} with $m=5$ demonstrates that $$\{1,4,6,9\}, \{2,3,7,8\}$$
form a $(10,2,4,2,8)$-DPDF and $(10,2,4,4,0)$-EPDF in $\mathbb{Z}_{10}$.
\item[(ii)] Applying Theorem \ref{directcyclic} with $m=9$ demonstrates that $$\{1,8,10,17\}, \{2,7,11,16\}, \{3,6,12,15\}, \{4,5,13,14\}$$
form a $(18,4,4,2,16)$-DPDF and $(18,4,4,12,0)$-EPDF in $\mathbb{Z}_{18}$.
\end{itemize}
\end{example}

Next, we present a construction in $\mathbb{Z}_{2p}$ where $p$ is a prime congruent to $1 \mod 4$.  It uses the fact that the non-zero squares in $GF(p)$ form a PDS when $p \equiv 1 \mod 4$ (\cite{Ma}). 

\begin{theorem}\label{DDStheorem}
Let $p$ be a prime congruent to $1 \mod 4$. 
Let $G$ be the additive group $\mathbb{Z}_{2p}$ and let $H=\{0,p\} \leq G$. \\ 
Define subsets $A_0, A_1$ of $\mathbb{Z}_{2p}$ as follows:
\begin{itemize}
    \item $A_0=\{s,s+p \in \mathbb{Z}_{2p}: s \mbox{ is a non-zero quadratic residue modulo $p$} \}$
    \item $A_1=\{t,t+p \in \mathbb{Z}_{2p}: t \mbox{ is a quadratic non-residue modulo $p$} \}$.
    \end{itemize}
Note $|A_0|=|A_1|=p-1$ and $\mathbb{Z}_{2p} \setminus \{A_0 \cup A_1\}=H$.  Then
\begin{itemize}
\item[(i)]
$$ \Delta(A_0)=\left(\frac{p-5}{2}\right) A_0 \cup \left(\frac{p-1}{2} \right) A_1 \cup (p-1) \{p\}.$$
\item[(ii)] 
$$\Delta(A_1)= \left(\frac{p-5}{2} \right) A_1 \cup \left(\frac{p-1}{2} \right) A_0 \cup (p-1) \{p\}.$$
\item[(iii)] $\{A_0,A_1\}$ forms a $(2p,2,p-1,p-3,2p-2)$-DPDF and a $(2p,2,p-1,p-1,0)$-EPDF in $\mathbb{Z}_{2p}$.
\end{itemize}
\end{theorem}
\begin{proof}
Let $Q_{2p}=\{a_1,\ldots,a_{\frac{p-1}{2}}\}$ be the quadratic residues modulo $p$, viewed as elements of $\mathbb{Z}_{2p}$; similarly, let $N_{2p}=\{b_1,\ldots,b_{\frac{p-1}{2}}\}$ be the quadratic non-residues modulo $p$, viewed as elements of $\mathbb{Z}_{2p}$.  For later convenience, we order the elements of $Q_{2p}$ in increasing order when viewed as integers, i.e. $0<a_1< a_2 < \cdots a_{\frac{p-1}{2}}<p$ as integers.  Then in $\mathbb{Z}_{2p}$, $$A_0=\{a_1, \ldots, a_{\frac{p-1}{2}}\} \cup \{ p+a_1, p+a_2, \ldots, p+a_{\frac{p-1}{2}}\}= Q_{2p} \cup (p+ Q_{2p})$$ 
where $p+Q_{2p}=\{p+x: x \in Q_{2p} \}$.

Consider $\Delta(A_0)$.  Clearly
$$ \Delta(A_0)= \Delta(Q_{2p}) \cup \Delta(Q_{2p}, p+Q_{2p}) \cup \Delta(p+Q_{2p}, Q_{2p}) \cup \Delta(p+Q_{2p}) $$
The multiset of internal differences of a set is unchanged by translation of the set, so $\Delta(Q_{2p})=\Delta(p+Q_{2p})$.  Since in $\mathbb{Z}_{2p}$
$$ a_i-(p+a_j)=-p+(a_i-a_j)=p+(a_i-a_j) \mbox{ and }(p+a_i)-a_j=p+(a_i-a_j)$$
we have that $\Delta(Q_{2p}, p+Q_{2p})=\{a_i-(p+a_j): a_i,a_j \in Q_{2p}\}=(p+\Delta(Q_{2p})+(\frac{p-1}{2}) \{p\}$ (in this multiset, unlike in $\Delta(Q_{2p})$, we have a contribution from terms with indices $i = j$).  A similar argument shows that $\Delta(p+Q_{2p},Q_{2p})=(p+\Delta(Q_{2p}))+(\frac{p-1}{2}) \{p\}$.  
Hence, to determine $\Delta(A_0)$, it suffices to determine $\Delta(Q_{2p})$.  

Denote by $Q_p$ the non-zero quadratic residues modulo $p$ viewed as elements of $\mathbb{Z}_p$.  It is well-known (see \cite{Ma}) that, as a subset of $\mathbb{Z}_p$, $Q_p$ is a $(p,\frac{p-1}{2},\frac{p-5}{4}, \frac{p-1}{4})$-PDS (where the element $\{0\}$ also occurs as an internal difference with frequency $\frac{p-1}{2}$).  Due to the order imposed on the elements of $Q_{2p}$, it is clear that the elements $D_>^{2p}=\{a_i-a_j: a_i>a_j\}$ of $\Delta(Q_{2p})$ will be precisely the same integers as the corresponding elements $D_>^p=\{a_i-a_j: a_i>a_j\}$ in $\Delta(Q_p)$, while the elements $D_<^{2p}=\{a_i-a_j: a_i<a_j\}$ of $\Delta(Q_{2p})$ will be $p+D_<^p$ where $D_<^p=\{a_i-a_j: a_i<a_j\}$ in $\Delta(Q_p)$.  Combining the multisets $\Delta(Q_{2p})$ and $\Delta(Q_{2p},p+Q_{2p})$ therefore yields $\frac{p-5}{4}$ copies of $Q_p \cup (p+Q_{2p})=A_0$, $\frac{p-1}{4}$ copies of $N_{2p} \cup (p+N_{2p})=A_1$, and $\frac{p-1}{2}$ copies of $\{p\}$.  Combining all multisets which make up $\Delta(A_0)$ yield the result.\\
A precisely analogous argument holds for the quadratic non-residues $N_{2p}$ to yield part (ii) (since $N_p$, the set of quadratic non-residues modulo $p$ viewed as elements of $\mathbb{Z}_{p}$, is also a $(p,\frac{p-1}{2},\frac{p-5}{4}, \frac{p-1}{4})$-PDS in $\mathbb{Z}_p$).  Finally, combining (i) and (ii) shows that the internal differences of $\{A_0,A_1\}$ comprise $\frac{p-5}{2}+\frac{p-1}{2}=p-3$ copies of each element of $A_0 \cup A_1=G \setminus H$, and $2(p-1)$ copies of the non-zero elements of $H$, hence form a DPDF with the stated parameters.  Since $G \setminus H$ is a $(2p,2p-2,2p-4,2p-2)$-PDS, the EPDF result follows.
\end{proof}

\begin{example}\label{DDSexample}
\begin{itemize}
    \item[(i)] For $p=5$, $A_0=\{1,4,6,9\}$ and $A_1=\{2,3,7,8\}$ in $\mathbb{Z}_{10}$ form a $(10,2,4,2,8)$-DPDF and a $(10,2,4,4,0)$-EPDF.
    \item[(ii)] For $p=13$, 
    $$A_0=\{1,3,4,9,10,12,14,16,17,22,23,25\}$$ 
    and
    $$A_1=\{2,5,6,7,8,11,15,18,19,20,21,24\}$$
    form a $(26,2,12,10,24)$-DPDF and a $(26,2,12,12,0)$-EPDF.
\end{itemize}
\end{example}

Observe that Example \ref{DDSexample}(i) is the same DPDF/EPDF obtained in Example \ref{directcyclicex}.  This is because, in $GF(5)$, the set of non-zero squares is $\{1,-1\}$.

In \cite{LeuMaTay}, a characterization is given for nontrivial reversible DDSs in cyclic groups (reversible means that $D=D^{-1}$ for the DDS $D$).  The result shows that (up to complementation and equivalence) there are only two possibilities for such a DDS $D$ and group $G$.  The first possibility is that $G=\mathbb{Z}_{2p}$ where $p$ is an odd prime with $p \equiv 1 \mod 4$ and $D$ is precisely $A_0 \cup \{0\}$ from Theorem \ref{DDStheorem}; $D$ is a $(p,2,p,p-1,\frac{p-1}{2})$-DDS in $G$ relative to $H=\{0,p\}$.  Our proof of Theorem \ref{DDStheorem} demonstrates directly how the partial difference set of quadratic residues mod $p$ gives the required properties for this DDS (whereas the proof in \cite{LeuMaTay} follows from a structural characterization using Sylow subgroups, combined with parameter restrictions from \cite{Ma}, so is not constructive).

Finally, we present an infinite family of DPDF/EPDFs constructed via coset partitioning.
\begin{theorem}\label{cosetpartition}
Let $G=\mathbb{Z}_{12d+4}$, where $d \in \mathbb{N}$, and define the subgroups $H = \{0,3d+1,6d+2,9d+3\} \cong \mathbb{Z}_4$ and $K=\{0,6d+2\} \cong \mathbb{Z}_2$. \\
Let $\{a_1+H, \ldots, a_{3d}+H \}$ be the cosets of $H$ in $G$, other than $H$ itself.\\
Define $S^{\prime}$ as follows:
\begin{itemize}
    \item partition each $a_i+H$ ($1 \leq i \leq d$) into two non-trival cosets of $K$, namely $a_i+K$ and $b_i+K$ where $b_i=a_i+(3d+1)$;
    \item partition each $a_j+H$ ($d+1 \leq j \leq 3d$) into two subsets $B_j$ and $C_j$, each of the form $\{d_k,d_l\}$ where $d_l-d_k = 3d+1$ (i.e. $\{a_j,a_j+(3d+1)\}$ and $\{a_j+(6d+2), a_j+(9d+3)\}$ or $\{a_j+(9d+3),a_j\}$ and $\{a_j+(6d+2), a_j+(9d+3)\}$);
    \item take $S^{\prime}$ to be this collection of $6d$ sets, each of size $2$.
\end{itemize}    
Then $S^{\prime}$ is a $(12d+4,6d,2,0,4d)$-DPDF and a $(12d+4,6d,2,12d-4,8d)$-EPDF.
\end{theorem}

\begin{proof}
It is clear that the sets of $S^{\prime}$ partition $G \setminus H$.
Since $\Delta(K) = 2(K\setminus\{0\})$ and the multiset of internal differences is unchanged by translation, for any $g+K$ ($g \in G$), the multiset $\Delta(g+K) = 2(K \setminus\{0\})$. So by partitioning each $a_i+H$ ($1 \leq i \leq d$) into $a_i+K$ and $b_i+K$, then computing $\Delta(a_i+K)$ and $\Delta(b_i+K)$, we obtain a collection of multisets comprising $4d$ copies of $\{6d+2\}$ in total. For each set of the form $B_j$ or $C_j$, we have $\Delta(B_j) = \Delta(C_j)= \{3d+1,9d+3\} = H \setminus K$. By partitioning each $a_j+H$ ($d+1 \leq j \leq 3d$) into $B_j$ and $C_j$, and for each computing $\Delta(B_j)$ and $\Delta(C_j)$, we obtain $4d$ copies of $H \setminus K$. Hence ${\rm Int}(S^{\prime})$ comprises $4d$ copies of $H \setminus \{0\}$ and zero copies of $G \setminus H$, so $S^{\prime}$ is a DPDF with the stated parameters.  Since $G$ is a $(12d+4,12d,12d-4,12d)$-PDS, $S^{\prime}$ is also an EPDF with the stated parameters.
\end{proof}

For a given $d$, Theorem \ref{cosetpartition} yields several (equally valid) DPDFs/EPDFs depending on the choices made for the sets. 
\begin{example}
Taking $d=1$ in Theorem \ref{cosetpartition} will produce a $(16,6,2,0,4)$-DPDF and $(16,6,2,14,8)$-EPDF which partition $G \setminus H$ where $G=\mathbb{Z}_{16}$ and $H=\{0,4,8,12\}$.
\begin{itemize}
    \item[(i)] From Example \ref{DPDFexample}(ii), one example is $\{1,9\},\{5,13\},\{14,2\},\{6,10\},\{15,3\},\{7,11\} $.
    \item[(ii)] A different example is $\{3,11\}, \{7,15\}, \{13,1\}, \{5,9\}, \{2,6\}, \{10,14\}$.
\end{itemize}
\end{example}

\section{RDS-based constructions for DPDFs/EPDFs}

In this section, we will show how relative difference sets can naturally be used to construct DPDFs. 

Relative difference sets were first introduced by Bose in \cite{Bos}, though they were not named as such; he presented his result as the ``affine analogue" of Singer's Theorem on difference sets.  The name and concept of RDS were formally introduced by Elliott and Butson in \cite{EllBut}.  The original construction of Bose for an RDS has parameters $(q+1, q-1, q, 1)$, and is couched in terms of finite geometry; a formulation in terms of finite fields is given in \cite{GanSpe}.  A more general result with parameters  $(\frac{q^r-1}{q-1}, q-1, q^{r-1}, q^{r-2})$ has been proved in various other ways, including via linear recurring sequences \cite{EllBut} and (in a particularly clear exposition) via linear functionals \cite{Leo}.

\subsection{Extending the Bose RDS construction to DPDFs/EPDFs}

Our first result demonstrates how Bose's original construction elegantly extends to a construction of a DPDF.  Each of the component sets is a Bose RDS.  We present it first in finite field terminology then outline the finite geometry viewpoint.

\begin{theorem}\label{Bose}
Let $q$ be a prime power and let $\alpha$ be a primitive element of $GF(q^2)$ with primitive polynomial $f$ over $GF(q)$. 
For each $\alpha^i \in GF(q^2)$ ($0 \leq i \leq q^2-2$), there exist $a_i, b_i \in GF(q^2)$ such that $\alpha^i=a_i+b_i \alpha$.\\

\begin{itemize}
    \item[(i)] For each $c \in GF(q)^*$, let $$S_c:=\{ \alpha^i \in GF(q^2)^*: \alpha^i=a_i+ c \alpha\}.$$
Then the family  $\{S_c\}_{c \in GF(q)^*}$ is a multiplicative $(q^2-1, q-1, q, q-1,0)$-DPDF and a multiplicative $(q^2-1,q-1,q,(q-1)(q-2), q^2-q)$-EPDF in $GF(q^2)^*$.
    \item[(ii)] For each $c \in GF(q)^*$, let $$S_c^{\prime}:=\{ i: \alpha^i \in S_c, 0 \leq i \leq q^2-2 \} \subseteq \mathbb{Z}_{q^2-1}.$$
Then the family $\{  S_c^{\prime} \}_{c \in GF(q)^*}$ is an additive $(q^2-1, q-1, q, q-1,0)$-DPDF and an additive $(q^2-1,q-1,q, (q-1)(q-2), q^2-q)$-EPDF in $\mathbb{Z}_{q^2-1}$.
    \end{itemize}
\end{theorem}
\begin{proof}
Let $c \in GF(q)^*$.  To construct the set $S_c$, we first form the multiplicative cosets of $GF(q)^*$ in $GF(q^2)^*$, and express their elements in the form $a+b \alpha$ ($a,b \in GF(q)$) via the primitive polynomial. 
There are $\frac{q^2-1}{q-1}=q+1$ cosets, each of the form $C_i$ where $C_0=\langle \alpha^{q+1} \rangle \cong GF(q)^*$ and $C_i=\alpha^i C_0$ ($0 \leq i \leq q$) (each of size $q-1$).  Each coset has the form $$C_i=\{ t \alpha^i: t \in GF(q)^*\}=\{ta_i+ tb_i \alpha: t \in GF(q)^*\}.$$
Observe that, in each $C_i$ other than $C_0$, there is a unique element  whose coefficient of $\alpha$ is $c$, namely $cb_i^{-1}(a_i+b_i\alpha)$.  Hence for each $c \in GF(q)^*$, we have
$$ S_c=\{cb_1^{-1}(a_1+b_1\alpha), cb_2^{-1}(a_2+b_2\alpha), \ldots, cb_q^{-1}(a_q+b_q\alpha)\}.$$
As $c$ runs through $GF(q)^*$, the $q-1$ sets $\{ S_c \}$ partition the elements of $GF(q^2)^* \setminus C_0$.

Now, each $S_c$ is a (multiplicative) $(q+1,q-1,q,1)$-RDS in $GF(q^2)^*$ with respect to the multiplicative subgroup $\langle \alpha^{q+1} \rangle$.  In fact, any one of these is a Bose RDS.  To see that no element of $C_0$ arises as a (multiplicative) difference, observe that each element of $S_c$ is in a distinct coset of $C_0$.  It can be shown by direct calculation in the finite field that every element of $GF(q^2) \setminus C_0$ arises precisely once as a difference, by considering the elements of $\Delta(S_c)$ directly (details are left to the reader).  Hence the family $\{ S_c\}_{c \in GF(q)^*}$ form a (multiplicative) $(q^2-1, q-1, q, q,0)$-DPDF and $(q^2-1,q-1,q, (q-1)(q-2), q^2-q)$-EPDF.

Finally, take the set of powers of $\alpha$ to convert each $S_c$ to a set $S_c^{\prime}=\{ i: \alpha^i \in S_c, 0 \leq i \leq q^2-2 \} \subseteq \mathbb{Z}_{q^2-1}$.  It is clear that this yields a collection $\{ S_c^{\prime} \}_{c \in GF(q)^*}$ of additive $(q+1,q-1,q,1)$-RDSs which form an additive DPDF and EPDF in $\mathbb{Z}_{q^2-1}$, with the same parameters as in the multiplicative case.
\end{proof}

This has a natural interpretation in finite geometry.  Consider $a + b \alpha \in GF(q^2)$ as the point $(a,b)$ in the affine plane $AG(2,q)$.  The construction above may be viewed as taking sets which are lines in a given parallel class (in this case, the class with $y=c$).  A parallel class has $q$ lines, each with $q$ points, and the lines in the class partition the points of the affine plane.

From Bose's paper, \cite{Bos}, each line with $c \neq 0$ in the parallel class gives an RDS with the required parameters, and taking all $q-1$ such lines, we obtain the DPDF described. To see this directly, replace each point in $AG(2,q)$ by the corresponding power of $\alpha$ via the above identification.  Since $(0,0)$ does not correspond to a power of $\alpha$, the line with $c=0$ is missing a point: the differences between all remaining points on this line are all multiples of $q+1$.  This corresponds to our omitted subgroup. For any other line, the differences will be precisely one occurrence of each of the $q(q-1)$ elements of $\mathbb{Z}_{q^2-1}$ that are not multiples of $q+1$, corresponding to our RDS. 

We note that the generalized RDS construction with parameters $(\frac{q^r-1}{q-1}, q-1, q^{r-1}, q^{r-2})$ cannot be used in this way to create a DPDF.

\begin{example}
Consider $GF(25)$ and let $\alpha$ be a primitive element, with primitive polynomial $x^2+x+2$ over $GF(5)$.  We have
$$ \alpha, \alpha^2=3+4 \alpha, \alpha^3=2+4 \alpha, \alpha^4=2+3 \alpha, \alpha^5=4+4\alpha, \alpha^6=2$$
Here $C_0=\{1=\alpha^0, 2=\alpha^6, 3=\alpha^{18}, 4=\alpha^{12}\}$, $C_1=\{ \alpha, 2\alpha, 3\alpha, 4\alpha \}$, $C_2=\{3+4\alpha, 1+3 \alpha, 4+2 \alpha, 2+ \alpha \}$, $C_3=\{ 2+ 4 \alpha, 4+3 \alpha, 1+ 2 \alpha, 3+\alpha \}$, $C_4=\{2+3 \alpha, 4+\alpha, 1+4 \alpha, 3+2 \alpha\}$, $C_5=\{4+4\alpha, 3+3 \alpha, 2+2 \alpha, 1+\alpha\}$.

Hence $S_1=\{\alpha^1, \alpha^{14}, \alpha^{15}, \alpha^{10}, \alpha^{17} \}$.  Taking the powers of $\alpha$, we obtain $S_1^{\prime}=\{1,14,15,10,17\}$: it can be checked that its internal differences comprise $1$ copy each of $\mathbb{Z}_{24} \setminus \{0,6,12,18\}$ and no copies of $\{6,12,18\}$.  \\
The other sets are obtained similarly:
\begin{itemize}
    \item $S_2=\{\alpha^7, \alpha^{20}, \alpha^{21}, \alpha^{16}, \alpha^{23}\} \mbox{ and }S_2^{\prime}=\{7,20,21,16,23\}$ 
 \item $S_3=\{\alpha^{19}, \alpha^{8}, \alpha^{9}, \alpha^{4}, \alpha^{11}\} \mbox{ and }S_3^{\prime}=\{19,8,9,4,11\}$
\item $S_4=\{\alpha^{13}, \alpha^{2}, \alpha^{3}, \alpha^{22}, \alpha^{5}\} \mbox{ and }S_4^{\prime}=\{13,2,3,22,5\}$.
\end{itemize}
Each $S_c$ is a $(6,4,5,1)$-RDS and $\{S_1,S_2,S_3,S_4\}$ is a $(24,4,5,4,0)$-DPDF and a $(24,4,5,12,20)$-EPDF.
\end{example}

\subsection{A general RDS construction for DPDFs/EPDFs}

In this section, we present an important general approach to constructing DPDFs and EPDFs using RDSs.  When dealing with groups which are not necessarily abelian, we will use multiplicative notation.

\begin{proposition}\label{RDSfirstprop}
Let $G$ be a group of order $mn$ and let $H$ be a (not necessarily normal) subgroup of $G$ of order $n$.
If $T=\{D_1, \ldots, D_s\}$ is a family of disjoint $k$-subsets of $G$ such that 
\begin{itemize}
    \item[(i)] each $D_i$ ($1 \leq i \leq s)$ is an $(m,n,k,\lambda)$-RDS in $G$ relative to $H$;
    \item[(ii)] $T$ partitions $G \setminus H$;
\end{itemize}
then $T$ is a $(mn,s,k,s\lambda,0)$-DPDF and an $(mn,s,k,mn-2n-s\lambda,mn-n)$-EPDF.
\end{proposition}
\begin{proof}
The multiset ${\rm Int}(T)$ comprises $s \lambda$ copies of each element of $G \setminus H$, and $0$ copies of $H$, and $T$ partitions $G \setminus H$, hence $T$ is an $(mn,s,k,s\lambda,0)$-DPDF. By Lemma \ref{PDSlemma}, $G \setminus H$ is an $(mn,mn-n,mn-2n,mn-n)$-PDS and so by Theorem \ref{PDSpartition}, $T$ is an $(mn,s,k,mn-2n-s\lambda, mn-n)$-EPDF.
\end{proof}

One natural method of producing a collection of disjoint sets with similar properties is to take an original set then form a collection of its translates by suitable group elements.  The next result, based on a lemma in \cite{GanSpe}, indicates what ``suitable" means in this context.

\begin{lemma}\label{RDSdisjoint}
Let $G$ be a group, $H$ a (not necessarily normal) subgroup of $G$ and let $D$ be an $(m,n,k,\lambda)$-RDS relative to $H$.  \\
Let $g_1 \neq g_2 \in G$.  The translates $g_1 D$ and $g_2 D$ of $D$ are disjoint if and only if $g_2^{-1} g_1 \in H$.
\end{lemma}
\begin{proof}
Suppose there is an element in $|g_1 D \cap g_2 D|$, i.e. $g_1d_1=g_2d_2$ for some $d_1,d_2 \in D$.  Then $g_2^{-1} g_1=d_2 d_1^{-1}$.  If $g_2^{-1}g_1 \in H$ then, since there is no non-identity element of $H$ of the form $d_2 d_1^{-1}$, we must have $g_2^{-1}g_1=1$ i.e. $g_1=g_2$, a contradiction.  So in this case $g_1 D$ and $g_2 D$ are disjoint.  Otherwise $g_2^{-1}g_1 \in G \setminus H$; this element has $\lambda>0$ representations in the multiset $\Delta(D)$.  So there is at least one pair $(d_1,d_2) \in D \times D$  such that $d_1 d_2^{-1}=g_2^{-1}g_1$.  Hence $g_2 d_1=g_1 d_2$ and so $|g_1 D \cap g_2 D|>0$.
\end{proof}

In order for a set of translates of $D$ to partition $G \setminus H$, we require $k|mn-n$.  Lemma \ref{RDSdisjoint} suggests translating by the elements of $H$; in this case we require $n+kn=mn$, i.e. $k=m-1$.  Combining this with the RDS relation $k(k-1)=(mn-n)\lambda$ (see comment following Definition \ref{RDSdefinition}), we have $(m-1)(m-2)=(m-1)n\lambda$, which implies $m=n \lambda+2$, hence $\lambda=\frac{m-2}{n}$. 

We present an explicit example of a DPDF/EPDF satisfying Proposition \ref{RDSfirstprop}, formed from an RDS relative to a subgroup $H$, translated by the elements of $H$.  It takes as its main ingredient the non-abelian RDS from \cite{CheLi} mentioned in a previous section.  This RDS is notable as being the first example of an RDS in a finite simple group with a non-trivial forbidden subgroup. The following DPDF/EPDF construction is noteworthy since it is the first known non-abelian example of a proper DPDF or EPDF.  Since the group is non-abelian, the result is written in multiplicative notation.

\begin{proposition}
Let $G$ be the alternating group $A_5$ acting on $\{1,2,3,4,5\}$ and let $\alpha=(25)(34)$.\\ Let $R$ be the following set of elements of $G$:
$$ (13542), (154), (14)(23), (13254), (12543), (15324), (245), (132),$$
$$(152), (13)(45), (12435), (235), (15234), (15342), (125), (14)(25),$$
$$(13)(25), (123), (14325), (23)(45), (14235), (253), (254), (13452),$$
$$ (12453), (145), (12)(35), (14523), (15)(24).$$

Let $R^{\prime}=\alpha R$:
$$ (13245), (15243), (14253), (13)(24), (124), (15423), (354), (13425),$$
$$ (15)(34), (13524), (12)(45), (345), (153), (15)(23), (12)(34), (143),$$
$$(134), (12534), (142), (24)(35), (14532), (234), (243), (135),$$
$$(12354), (14352), (12345), (14)(35), (15432).$$
Then
\begin{itemize}
\item[(i)] $R$ is a $(30,2,29,14)$-RDS in $G$ relative to the subgroup $H=\langle \alpha \rangle \cong \mathbb{Z}_2$.
\item[(ii)] $A=\{ R, R^{\prime}\}$ is a $(60,2,29,28,0)$-DPDF and $(60,2,29,28,58)$-EPDF which partitions $G \setminus H$.
\end{itemize}
\end{proposition}
\begin{proof}
For (i), the proof that $R$ is a $(30,2,29,14)$-RDS is the content of the paper \cite{CheLi}.  It is obtained using structural properties of Cayley graphs.\\
For (ii), we show it satisfies Proposition \ref{RDSfirstprop}.  Here $m=30$, $n=2$, $s=2$, and $R^{\prime}=\alpha R$ (the translate of $R$ by $\alpha$). By inspection, neither element of $H=\{id, \alpha\}$ is contained in $R$ nor $R^{\prime}$, and these two sets partition $G \setminus H$.  By (i), $R$ is an  $(30,2,29,14)$-RDS relative to $H$.  To see that the same is true of $R^{\prime}$, observe that  $\Delta(R^{\prime})$ is the multiset of all elements of the form $(\alpha r_1)(\alpha r_2)^{-1}=\alpha (r_1 r_2^{-1}) \alpha$ ($r_1,r_2 \in R$), i.e. the multiset $\{\alpha x \alpha: x \in \Delta(R)\}$.  Now, $\Delta(R)$ comprises $14$ copies of $G \setminus H$ and $0$ copies of $H$. We have $\alpha G \alpha=G$ and $\alpha H \alpha=\alpha \{1,\alpha \} \alpha =H$, so $\alpha (G \setminus H) \alpha=G \setminus H$ and hence $\Delta(R^{\prime})=\Delta(R)$. Hence by Proposition \ref{RDSfirstprop}, $A$ is a $(60,2,29,28,0)$-DPDF and a $(60,2,29,28,0)$-DPDF and a $(60,2,29,28,58)$-EDPF in $A_5$.
\end{proof}

In general, it is not feasible to perform explicit verification of the properties required for Proposition \ref{RDSfirstprop}: we will therefore establish results which guarantee that large classes of structures satisfy the requirements of Proposition \ref{RDSfirstprop}.  Henceforth we will assume that $H$ is a normal subgroup of $G$.

\begin{lemma}\label{RDSbasic}
Let $G$ be a group, $H$ a normal subgroup of $G$ and let be $D$ an $(m,n,k,\lambda)$-RDS relative to $H$. 
\begin{itemize}
\item[(i)] $\Delta(D)=\Delta(gD)$ for any $g \in G$. \\
 In particular, for any $g \in G$, any translate $gD$ is a $(m,n,k,\lambda)$-RDS relative to $H$.
\item[(ii)] $D$ cannot contain more than one representative from any coset $gH$ ($g \in H$).  In particular, $k \leq m$.
\end{itemize}
\end{lemma}
\begin{proof}
 For (i), $gd_1 (g d_2)^{-1}=g(d_1 d_2^{-1})g^{-1}$.  By definition, the multiset $\Delta(D)$ comprises $\lambda$ copies of $G \setminus H$ and $0$ copies of $H$. Since $H$ is a normal subgroup, for any $g \in G$ we have $gHg^{-1}=H$, and since $gGg^{-1}=G$ we also have that $g(G \setminus H)g^{-1}=G \setminus H$. So the multiset $\Delta(gD)$ also has $\lambda$ copies of $G \setminus H$ and $0$ copies of $H$.\\
(ii) Suppose $D$ contains elements $d_1 \neq d_2 \in gH$, say $d_1=gh_1$ and $d_2=gh_2$. Then $d_1 d_2^{-1}=gh_1(gh_2)^{-1}=g(h_1 h_2^{-1})g^{-1} \in H$ since $H$ is normal in $G$.  Since $D$ is an RDS relative to $H$, we must have $d_1 d_2^{-1}=1$, i.e. $d_1=d_2$, a contradiction.
\end{proof}

\begin{theorem}\label{RDSgeneral}
Let $G$ be a group of order $mn$, let $H$ be a normal subgroup of $G$ of order $n$, and suppose there exists an $(m,n,m-1, \frac{m-2}{n})$-RDS $R$ in $G$ relative to $H$. \\
Then there exists a $(mn,n,m-1,m-2,0)$-DPDF and an $(mn,n,m-1,(m-2)(n-1),(m-1)n)$-EPDF which partitions $G \setminus H$.
\end{theorem}
\begin{proof}
Suppose we have an $(m,n,m-1,\frac{m-2}{n})$-RDS $R$.  Since $|R|=k=m-1=[G:H]-1$, and by Lemma \ref{RDSbasic} $R$ cannot contain a representative of more than one coset of $H$, there is precisely one coset of $H$ with no representative in $R$.  Without loss of generality, we may replace $R$ by a suitable translate $D:=gR$ ($g \in G$), so that the coset without a representative in $D$ is $H$ itself.  (This can be the trivial translation by the identity if $H \cap R=\emptyset$.)  By a previous result, any translate of $R$ is also an $(m,n,m-1,\frac{m-2}{n})$-RDS and has its elements in distinct cosets of $H$.  Hence $D$ is a $(m,n,m-1,\frac{m-2}{n})$-RDS comprising a representative of each coset of $H$ except $H$ itself.

Let $\mathcal{D_H}=\{hD: h \in H\}$; we shall show this is the desired DPDF/EPDF.
Since $D$ contains no element of $H$, any translate $hD$ with $h \in H$ must also have empty intersection with $H$ (if $h_1 \in (hD) \cap H$ then $h_1=hd$ for some $d \in D$, i.e. $d=h^{-1} h_1 \in H$, impossible).  By Lemma \ref{RDSdisjoint}, the sets in $\mathcal{D_H}$ are pairwise disjoint, i.e. their union comprises $kn$ distinct elements of $G$. Hence the sets of $\mathcal{D_H}$ partition the $mn-n=kn$ elements of $G \setminus H$.

By Lemma \ref{RDSbasic}, each set in $\mathcal{D_H}$ is an $(m,n,m-1,\frac{m-2}{n})$-RDS relative to $H$.  Finally, by Proposition \ref{RDSfirstprop} $\mathcal{D_H}$ is an $(mn,n,m-1,m-2,0)$-DPDF and an $(mn,n,m-1,(m-2)(n-1),(m-1)n)$-EPDF.
\end{proof}

\begin{example}
Let $G=\mathbb{Z}_8$ and consider the $(4,2,3,1)$-RDS $D=\{1,6,7\}$ relative to the subgroup $H=\{0,4\}$.  Note that the coset of $H$ not represented in $D$ is $H$ itself.  Then $\mathcal{D_H}=\{ \{1,6,7\}, \{5,2,3\} \}$ forms an $(8,2,3,2,0)$-DPDF and an $(8,2,3,2,6)$-EPDF.
\end{example}

\begin{remark}
Observe that the construction in Theorem \ref{Bose} extending the Bose approach uses a component RDS with parameters $(q+1,q-1,q,1)$, which satisfies the requirements of Theorem \ref{RDSgeneral}.  This is not a coincidence; we can view the Bose approach as an instance of Theorem \ref{RDSgeneral}, in the following way.

In the notation of Theorem \ref{Bose}, for each $\alpha^i$ ($1 \leq i \leq q$), we have $\alpha^i=a_i + b_i \alpha$ ($a_i,b_i \in GF(q)$).  The multiplicative coset $C_i$ has the form $C_i=\{ ta_i + tb_i \alpha: t \in GF(q)^*\}$, and so for any $c \in GF(q)^*$, the element of $C_i$ which lies in $S_c$ is $cb_i^{-1} \alpha^i$.  Hence we can write
$$ S_c=\{cb_1^{-1}\alpha, cb_2^{-1}\alpha^2, \ldots, cb_q^{-1} \alpha^q \}.$$
Taking discrete logs (viewed as elements of $\mathbb{Z}_{q^2-1}$), we have 
$$S_c^{\prime}=\log(S_c)=\log(c)+\{ \log(b_1^{-1}\alpha), \log(b_2^{-1}\alpha^2),\ldots, \log(b_q^{-1} \alpha^q) \}=\log(c)+S_1^{\prime}.$$
So the $(q^2-1,q-1,q,q-1,0)$-DPDF and $(q^2-1,q-1,q,(q-1)(q-2),q^2-q)$-EPDF in $\mathbb{Z}_{q^2-1}$ obtained in Theorem \ref{Bose} from the extension of the Bose approach may be viewed as $\mathcal{D_H}$ where $D=S_1^{\prime}$ and $H=\log(C_0) \cong \mathbb{Z}_{q-1}$, the excluded subgroup.
\end{remark}

\subsection{Applications of the general RDS construction}
An RDS with parameters $(n+1,n-1,n,1)$ (and $H$ normal in $G$) can always be used in the construction of Theorem \ref{RDSgeneral}.  An RDS with these parameters is said to be \emph{affine}; more detail about affine RDSs is given in \cite{GanSpe} and \cite{Pot}, including non-abelian examples.  It is conjectured that in the abelian case, $n$ must be a prime power.

The following existence result for a non-abelian affine RDS is from \cite{GanSpe}.
\begin{proposition}\label{Spence}
Let $n=p^r$ where $p$ is prime and $(G,+)$ is the cyclic group of integers modulo $p^{2r}-1$.  We define a new addition on the elements of $G$.  Let $q=p^h$ and suppose $hv=2r$ where $v$ is an integer all of whose prime factors divide $q-1$.  (If $q \equiv 3 \mod 4$, we also need $v \not\equiv 0 \mod 4)$.  Then, given $j$, let $r(j)$ denote the unique integer $i (\mod v)$ such that 
$$ q^i\equiv 1+j(q-1) \mod v(q-1).$$
Define a sum $\oplus$ on $G$ by
$$ i \oplus j \equiv iq^{r(j)}+j \mod q^v-1.$$
Then
\begin{itemize}
\item[(i)] $(G, \oplus)$ is a non-abelian group if $v>1$.
\item[(ii)] Let $n=p^r$.  If D is an $(n+1,n-1,n,1)$-RDS in $(G,+)$ relative to $(H,+)$, and the automorphism $x \mapsto px$ fixes $D$, and $(H, \oplus)$ is a normal subgroup of $(G, \oplus)$, then $D$ is a RDS in $(G,\oplus)$ relative to $(H,\oplus)$.
\item[(iii)] When $v=2$ and $p$ is odd, the conditions of (ii) are satisfied and so there exist examples of non-abelian $(n+1,n-1,n,1)$-RDSs.
\end{itemize}
\end{proposition}

By Theorem \ref{RDSgeneral}, since $(H,\oplus)$ is required to be a normal subgroup of $(G,\oplus)$ in Proposition \ref{Spence}, any RDS from the above construction guarantees the existence of non-abelian $(n^2-1,n-1,n,n-1,0)$-DPDFs and $(n^2-1,n-1,n,(n-1)(n-2),n(n-1))$-EPDFs.  

Next, we seek to identify RDS constructions which satisfy the requirements of Theorem \ref{RDSgeneral} but are not of affine type.  

We present a general construction for DPDFs and EPDFs using difference sets with additional properties.  This result extends and generalizes ideas from \cite{FanLeiSha} and \cite{GeFujMia}, which use such structures in cyclic groups to build optimal frequency-hopping sequences and difference systems of sets.

\begin{theorem}\label{DS_subgroup}
Let $n \in \mathbb{N}$ and let $m=(n-1)^2+1$.\\
Let $G$ be a group of order $mn$ and let $H$ be a normal subgroup of order $n$. \\
Suppose $D$ is a $(mn, (n-1)^2+n,n)$-difference set containing $H$. Let $R = D\backslash{H}$. Then
\begin{itemize}
\item[(i)] $R$ is an $(m,n,m-1, n-2)$-RDS relative to $H$; 
\item[(ii)] the sets 
$$\mathcal{R_H}=\{hR: h \in H\}$$ 
form an $(n[(n-1)^2+1],n,(n-1)^2,n(n-2),0)$-DPDF and an $(n[(n-1)^2+1],n,(n-1)^2,n(n-1)(n-2),n(n-1)^2)$-EPDF in $G$.
\end{itemize}
\end{theorem}

\begin{proof}
We show that $R$ is an RDS with the stated parameters.  Part (ii) then follows by observing $\frac{m-2}{n}=\frac{(n-1)^2-1}{n}=n-2$ and applying Theorem \ref{RDSgeneral}.

We have that
\begin{equation*}
    \Delta(D) = \Delta(R) \cup \Delta(R,H) \cup \Delta(H,R) \cup \Delta(H);
\end{equation*}
this multiset union comprises $n$ copies of each non-identity element of $G$.  Since $H$ is a subgroup of order $n$, the multiset $\Delta(H)$ comprises $n$ copies of each non-identity element of $H$ and no elements of $G \setminus H$.  Thus each of the $n(n-1)^2$ elements of $G\backslash{H}$ must occur precisely $n$ times across the multiset union $\Delta(R) \cup \Delta(R,H) \cup \Delta(H,R)$, and this accounts for all its $n^2(n-1)^2$ elements.

Non-identity elements of $H$ are obtained as differences $xy^{-1}$ in $\Delta(G)$ precisely when $x$ and $y$ are distinct and lie in the same coset of $H$. Since $\Delta(R) \cup \Delta(R,H) \cup \Delta(H,R)$ comprises $n$ copies of ${G}\backslash{H}$ and no copies of $H\backslash\{0\}$, it is clear that $R$ consists of at most one representative of each non-trivial coset $aH$ of $H$ ($R$ is disjoint from $H$ by construction). Since there are $(n-1)^2$ such cosets of $H$ and as $R$ has cardinality $(n-1)^2$, $R$ must consist of exactly one representative from each non-trivial coset of $H$. For any $g \in aH$, where $a \not\in H$, the multiset $\Delta(g,H) = aH$. From the structure of $R$, we see that the multiset $\Delta(R,H) = {G} \backslash H$. Analogously $\Delta(H,R) = {G} \backslash H$. Thus, $\Delta(R)$ comprises precisely $n-2$ copies of ${G}\backslash{H}$. So $R$ is an $((n-1)^2+1,n,(n-1)^2, n-2)$-RDS relative to subgroup $H$.
Part(ii) now follows by application of Theorem \ref{RDSgeneral}.
\end{proof}

\begin{example}
We present both a non-abelian and an abelian example of a $(40,13,4)$-DS containing a normal subgroup of order $4$ which can be used in Theorem \ref{DS_subgroup}.
\begin{itemize}
    \item[(i)] Let $G$ be the semi-direct product of $C_5$ and $C_8$ acting via $C_8/C_4=C_2$ (\cite{Dok}).\\
$G$ is a non-abelian group with presentation $\langle a,b: a^5=b^8=1, ba=a^4b \rangle$.
Its centre (a normal subgroup of $G$ by definition) is given by $H=\{1,b^2,b^4,b^6\} \cong C_4$.\\
From \cite{ColDin}, a $(40,13,4)$-DS is given by
$$ D= \{a^4,1,a,a^4b, a^2b,b^2,a^2b^2,a^3b^2,b^4,b^5,ab^5,b^6,a^3b^7\}.$$
Here $H \subseteq D$ and 
$$R=D \setminus H=\{a,a^4, a^2b, a^4b,a^2b^2,a^3b^2,b^5, ab^5,a^3b^7\}$$ is a $(10,4,9,2)$-RDS relative to $H$.  By Theorem \ref{DS_subgroup}, it yields a non-abelian $(40,4,9,8,0)$-DPDF and $(40,4,9,24,36)$-EPDF.
\item[(ii)] Let $G=\mathbb{Z}_{40}$.  A $(40,13,4)$-difference set (from \cite{ColDin}, arising from $PG(3,3)$) is given by 
$$D=\{0,6,7,8,10,11,14,19,20,23,25,30,32\}$$
and it contains the subgroup $H=10 \mathbb{Z}_4=\{0,10,20,30\}$.  Here $H \subseteq D$
 and $$R=D \setminus H=\{6,7,8,11,14,19,23,25,32\}$$
is a $(10,4,9,2)$-RDS in $\mathbb{Z}_{40}$.  By Theorem \ref{DS_subgroup}, it yields an abelian $(40,4,9,8,0)$-DPDF and $(40,4,9,24,36)$-EPDF.
\end{itemize}
\end{example}

Using a difference set result from \cite{GeFujMia}, we can  guarantee the existence of an infinite family of DPDFs and EPDFs in cyclic groups via Theorem \ref{DS_subgroup}.  In fact, \cite{GeFujMia} provides an explicit construction for an appropriate family of cyclic difference sets (which possess additional properties not required for our application) using finite geometry; we refer the reader to that paper for more details.

\begin{corollary}\label{GeFujMiaExample}
Let $n = 2^r + 1$ where $r$ is a positive integer and let $m=(n-1)^2+1$. \\
There exists an $(mn,n,(n-1)^2,n(n-2),0)$-DPDF and an $(mn,n,(n-1)^2,n(n-1)(n-2),n(n-1)^2)$-EPDF in $\mathbb{Z}_{mn}$.
\end{corollary}
\begin{proof}
It is proved in \cite{GeFujMia} that for $n=2^r+1$ there exists an $(n((n-1)^2 + 1), n + (n-1)^2, n)$-difference set $D$ in $\mathbb{Z}_{n((n-1)^2+1)}$ which contains  subgroup $H=((n-1)^2 + 1) \mathbb{Z}_n$.
The result follows by applying Theorem \ref{DS_subgroup}.
\end{proof}

\begin{example}
Applying Corollary \ref{GeFujMiaExample} using the difference set from $\cite{GeFujMia}$ with $n=5$ yields a cyclic $(85,21,5)$-DS given by
$$ D=\{0,1,2,4,7,8,14,16,17,23,27,28,32,34,43,46,51,54,56,64,68\}$$
which contains the subgroup $H=17 \mathbb{Z}_5=\{0,17,34,51,68\}$.\\
Let $R=D \setminus H=\{1,2,4,7,8,14,16,23,27,28,32,43,46,54,56,64\}$.\\
Then $\mathcal{R_H}=\{R,R+17,R+34,R+51,R+68\}$ is an $(85,5,16,15,0)$-DPDF and a $(85,5,16,60,80)$-EPDF in $\mathbb{Z}_{85}$.
\end{example}

We end this section by observing that not all DPDFs/EPDFs with the $\mu$-value of the DPDF equal to zero must have constituent sets which are RDSs; see for example the construction below.

\begin{example}
Let $\mathbb{Z}_{3m}$ where $m>3$ is an odd number, and let $H=\{0,m,2m\} \cong \mathbb{Z}_3$. \\
Let $S_i= \{i,3m-i\}$ ($1 \leq i \leq 3m-1)$.  It is straightforward to verify that $\Delta(S_i)=S_{2i}$ and that the family
$S^{\prime}=\{ S_i: 1 \leq i \leq \frac{3m-1}{2}, i \neq m\}$ form a $(3m,\frac{3m-3}{2},2,1,0)$-DPDF and a $(3m,\frac{3m-3}{2},2,3m-7,3m-3)$-EPDF.\\
When $m=5$, then the sets $\{1,14\}, \{2,13\}, \{3,12\}, \{4,11\}, \{6,9\},\{7,8\}$ form a $(15,6,2,1,0)$-DPDF and a $(15,6,2,8,12)$-EPDF in $\mathbb{Z}_{15}$.
\end{example}

\section{DPDFs which are not EDPFs and vice versa}

Although the main focus of the paper has been to consider structures which are simultaneously DPDFs and EPDFs, we end with some examples which show there exist DPDFs which are not EPDFs, and EPDFs which are not DPDFs. Our examples occur in cyclic groups.

\begin{proposition}
Let $n=2t+1$ be an odd number, $t>2$.  Then in $\mathbb{Z}_{2t+1}$, the family of sets $S^{\prime}$
$$ \{2,3\}, \{4,5\}, \ldots, \{2t-2,2t-1\}$$
is a $(2t+1,t-1,2,0,t-1)$-DPDF which is not an EPDF for $t>2$.
\end{proposition}
\begin{proof}
It is immediate to check for each set $S_i=\{2i,2i+1\}$ in $S^{\prime}$ ($1 \leq i \leq t-1$) the multiset $\Delta(S_i)=\{-1,+1\}$ while $S=G^* \setminus \{-1,+1\}$.  Hence ${\rm Int}(S^{\prime})$ comprises $t-1$ copies of $\{-1,+1\}$ and $0$ copies of $S=G^* \setminus \{-1,+1\}$. \\
The multiset ${\rm Int}(S)$ contains $\{-1,1\}$ $n-4$ times, $\{-2,2\}$ $n-5$ times and all other elements of $G^*$ $n-6$ times, hence $S$ is not a PDS and so $S$ is not an EPDF.
\end{proof}

The existence of an infinite family of EPDFs which are not DPDFs is still an open question (see next section).  The following examples in cyclic groups (\cite{Chris}, obtained via computation in GAP \cite{GAP}) show that there exist EPDFs which are not DPDFs:
\begin{example}
\begin{itemize}
    \item[(i)] $\{1,8\}, \{3,6\}$ is a $(9,2,2,0,2)$-EPDF in $\mathbb{Z}_9$ which is not a DPDF;
    \item[(ii)] $\{1,2,11,12\}, \{3,5,8,10\}$ is a $(13,2,4,2,4)$-EPDF in $\mathbb{Z}_{13}$ which is not a DPDF.
\end{itemize}
\end{example}

\section{Conclusions and further work}
This paper has demonstrated how the recently-introduced combinatorial structures of DPDFs and EPDFs can be constructed in groups other than the additive group of a finite field, using techniques other than the cyclotomic approach introduced in \cite{HucJoh}.  In particular, it has demonstrated the existence of infinite families and individual examples in non-cyclotomic abelian groups and non-abelian groups.  In general, our constructions partition the complement of a subgroup.

All constructions in this paper for $(v,s,k,\lambda_1,\mu_1)$-DPDFs which are $(v,s,k,\lambda_2,\mu_2)$-EPDFs have the property that at least one of the frequencies $\{\lambda_1,\lambda_2,\mu_1,\mu_2\}$ takes a zero value (even for constructions not using RDSs). This is in contrast to the cyclotomic constructions of \cite{HucJoh}, where there were numerous examples with all four frequencies non-zero.  This motivates the following:

\begin{openproblem}
Is it possible to obtain families of sets which are both $(v,s,k,\lambda_1,\mu_1)$-DPDFs and $(v,s,k,\lambda_2,\mu_2)$-EPDFs, not corresponding to cyclotomic constructions in the additive group of a finite field, such that $\{\lambda_1,\lambda_2,\mu_1,\mu_2\}$ are all non-zero? 
\end{openproblem}

To contextualize this Open Problem, we classify below the types of DPDF/EPDF construction partitioning $G \setminus H$ which we currently know, in terms of the nature of $\{\lambda_1,\lambda_2,\mu_1,\mu_2\}$:

\begin{theorem}\label{categories}
Let $S^{\prime}$ be both a $(v,s,k,\lambda_1,\mu_1)$-DPDF and a $(v,s,k,\lambda_2,\mu_2)$-EPDF which partitions $G \setminus H$, $H$ a normal subgroup of $G$.  Then
\begin{itemize}
\item[(i)] If $S^{\prime}$ consists of all non-trivial cosets of $H$ then $\lambda_1=0$ and $\mu_2=0$.
\item[(ii)] If every set in $S^{\prime}$ is a union of at least $2$ cosets of $H$, then $\lambda_1>0$, $\mu_1>0$ and $\mu_2=0$.
\item[(iii)] If the sets of $S^{\prime}$ are a subdivision of the nontrivial cosets of $H$ (i.e. formed by partitioning the cosets) then $\lambda_1=0$, $\mu_1>0$ and $\mu_2>0$.
\item[(iv)] If every set in $S^{\prime}$ has at most one representative from each coset of $H$, then $\mu_1=0$, $\lambda_1>0$ and $\mu_2>0$.
\end{itemize}
\end{theorem}

Examples of these types are as follows: for (i), see Theorem \ref{DPDFcosets}; for (ii), see Theorem \ref{directcyclic}, and for (iii) see Theorem \ref{cosetpartition}.  Illustrations of type (iv) include all RDS-based examples in Section 4.

Hence, any construction satisfying the Open Problem would need to lie outside the list of the categories given in Theorem \ref{categories} (as well as having $n>2$ and $gcd(sk,v-1)>1$ by Theorem \ref{GsetminusH}).  We note that these account for all examples of DPDFs/EPDFs partitioning the complement of a subgroup that the authors are currently aware of.

In the final section of this paper, we have provided a brief indication that there exist DPDFs which are not EPDFs, and vice versa.  It would be of interest to find more examples of such structures, in a variety of groups.  In particular:
\begin{openproblem}
Construct an infinite family of EPDFs which are not DPDFs.
\end{openproblem}

\subsection*{Acknowledgements}
We thank Maura Paterson for insights into finite geometry, and Jim Davis for discussions on relative difference sets.   We thank Chris Jefferson and Matthew McIlree for helpful computation. The research of Sophie Huczynska is supported by EPSRC grant EP/X021157/1.

\end{document}